\documentclass[12pt]{article}
\usepackage{amssymb}
\newcommand{\dd}{{\rm \kern 3pt I\kern-9pt d}}

\newcommand{\Abar}{{\backslash\kern-8pt A}}

\title{\large ON ERROR OPERATORS RELATED TO THE ARBITRARY FUNCTIONS PRINCIPLE}
\author{\sc Nicolas Bouleau\footnote{Ecole des Ponts, ParisTech, 28 rue des Saints P\`eres, 75007 Paris, France; 
e-mail : {\tt bouleau@enpc.fr}}}
\date{October 2006}
\newcommand{\CC}{\mbox{{\helv \i}{\rm\kern-5pt C}}}

\begin{document}
\maketitle

\noindent{\bf Abstract.} The error on a real quantity $Y$ due to the graduation of the measuring instrument may be asymptotically represented, when the graduation is regular and fines down, by a Dirichlet form on $\mathbb{R}$ whose square field operator does not depend on the probability law of $Y$ as soon as this law possesses a continuous density. This feature is related to the ``arbitrary functions principle" (Poincar\'e, Hopf). We give extensions of this property to $\mathbb{R}^d$ and to the Wiener space for some approximations of the Brownian motion. This gives  new approximations of the Ornstein-Uhlenbeck gradient. These results apply to the discretization of some stochastic differential equations encountered in mechanics. \\


\noindent{\Large\textsf{Introduction.}}\\

\noindent The approximation of a random variable $Y$ by an other one $Y_n$ yields most often a Dirichlet form.  The framework is general, cf. Bouleau [3]  whose results are recalled \S I.1 below.

Usually, when this Dirichlet form exists and does not vanish, the conditional law of $Y_n$ given $Y=y$ is not reduced to a Dirac mass, and the variance of this conditional law yields the square field operator $\Gamma$. On the other hand when the approximation is deterministic, i.e. when $Y_n$ is a function of $Y$ say $Y_n=\eta_n(Y)$, then {\it most often} the symmetric bias operator $\tilde{A}$ and the Dirichlet form vanish, cf. Bouleau [3] examples 2.1 to 2.9 and remark 5.

Nevertheless, there are cases where the conditional law of $Y_n$ given $Y$ is a Dirac mass, i.e. $Y_n$ is a deterministic function of $Y$, and where the approximation of $Y$ by $Y_n$ yields even so a non zero Dirichlet form on $L^2(\mathbb{P}_Y)$.

This phenomenon is interesting, insofar as randomness (here the Dirichlet form) is generated by a deterministic device. In its simplest form, the phenomenon appears precisely when a quantity is measured by a graduated instrument to the nearest graduation when looking at the asymptotic limits as the graduation fines down. 

The first part of this article is devoted to functional analytical tools that we need afterwards. We first recall the properties of
the bias operators and  the Dirichlet form associated with an approximation. Next we prove a Girsanov-type theorem for Dirichlet forms which has its own interest, i.e. an answer to the question of an absolutely continuous change of measure for Dirichlet forms. At last we recall some simple properties of Rajchman measures.

The second part is devoted to the case of a real or finite dimensional quantity measured with equidistant graduations. The mathematical argument here is basically {\it the arbitrary functions method} about which we give a short historical comment.

Several infinite dimensional extensions of the arbitrary functions principle are studied in the third part. The first one is about approximations of continuous martingales whose brackets are Rajchman measures. Then we consider the case of the Wiener space on which the preceding results may be improved and other asymptotic properties are obtained concerning the approximation of the Ornstein-Uhlenbeck gradient. Eventually we apply these results to the approximation of stochastic differential equations encountered in mechanics and solved by the Euler scheme.\\

\noindent{\Large\textsf{I. Functional analytical tools.}}\\

\noindent{\large\textsf{I.1. }}{\bf Approximation, Dirichlet forms and bias operators.}

Our study uses the theoretical framework concerning the bias operators and the Dirichlet form generated by an approximation proposed in Bouleau [3]. We recall here the definitions and main results for the convenience of the reader. Here, considered Dirichlet forms are always symmetric.

Let $Y$ be a random variable  defined on $(\Omega, \mathcal{A}, \mathbb{P})$ with values in a measurable space $(E,\mathcal{F})$ and let $Y_n$ be approximations also defined on $(\Omega, \mathcal{A}, \mathbb{P})$ with values in $(E,\mathcal{F})$. We consider an algebra $\mathcal{D}$ of bounded functions from $E$ into $\mathbb{R}$ or $\mathbb{C}$ containing the constants and dense in $L^2(E,\mathcal{F},\mathbb{P}_Y)$ and a sequence $\alpha_n$ of positive numbers. With $\mathcal{D}$ and $(\alpha_n)$ we consider the four following assumptions defining the four bias operators 
$$
(\mbox{H}1)\qquad\left\{\begin{array}{l}
\forall \varphi\in{\cal D}, \mbox{ there exists } \overline{A}[\varphi]\in L^2(E,{\cal F},\mathbb{P}_Y)\quad s.t. \quad\forall \chi\in{\cal D}\\
\lim_{n\rightarrow\infty} \alpha_n\mathbb{E}[(\varphi(Y_n)-\varphi(Y))\chi(Y)]=\mathbb{E}_Y[\overline{A}[\varphi]\chi].
\end{array}\right.
$$
$$
(\mbox{H}2)\qquad\left\{\begin{array}{l}
\forall \varphi\in{\cal D}, \mbox{ there exists } \underline{A}[\varphi]\in L^2(E,{\cal F},\mathbb{P}_Y)\quad s.t. \quad\forall \chi\in{\cal D}\\
\lim_{n\rightarrow\infty} \alpha_n\mathbb{E}[(\varphi(Y)-\varphi(Y_n))\chi(Y_n)]=\mathbb{E}_Y[\underline{A}[\varphi]\chi].
\end{array}\right.
$$
$$
(\mbox{H}3)\quad\left\{\begin{array}{l}
\forall \varphi\in{\cal D}, \mbox{ there exists } \widetilde{A}[\varphi]\in L^2(E,{\cal F},\mathbb{P}_Y)\quad s.t. \quad\forall \chi\in{\cal D}\\
\lim_{n\rightarrow\infty} \alpha_n\mathbb{E}[(\varphi(Y_n)-\varphi(Y))(\chi(Y_n)-\chi(Y))]=-2\mathbb{E}_Y[\widetilde{A}[\varphi]\chi].
\end{array}\right.
$$
$$
(\mbox{H}4)\quad\left\{\begin{array}{l}
\forall \varphi\in{\cal D}, \mbox{ there exists } \Abar[\varphi]\in L^2(E,{\cal F},\mathbb{P}_Y)\quad s.t. \quad\forall \chi\in{\cal D}\\
\lim_{n\rightarrow\infty} \alpha_n\mathbb{E}[(\varphi(Y_n)-\varphi(Y))(\chi(Y_n)+\chi(Y))]=2\mathbb{E}_Y[\Abar[\varphi]\chi].
\end{array}\right.
$$
We first note that as soon as two of hypotheses (H1) (H2) (H3) (H4) are fulfilled (with
 the same algebra ${\cal D}$ and the same sequence $\alpha_n$), the other two follow thanks to the relations
$$\widetilde{A}=\frac{\overline{A}+\underline{A}}{2}\quad\quad\Abar=\frac{\overline{A}-\underline{A}}{2}.$$
When defined, the operator $\overline{A}$ which considers the asymptotic error from the point of view of the limit model, will be called {\it the 
theoretical bias operator}.

The operator $\underline{A}$ which considers the asymptotic error from the point of view of the approximating model will be called {\it the 
practical bias operator}.

Because of the property
$$\langle\widetilde{A}[\varphi],\chi\rangle _{L^2(\mathbb{P}_Y)}=\langle\varphi,\widetilde{A}[\chi]\rangle _{L^2(\mathbb{P}_Y)}$$
the operator $\widetilde{A}$ will be called {\it the symmetric bias operator}.

The operator $\Abar$ which is often (see theorem 2 below) a first order operator will be called {\it the singular bias operator}.\\

\noindent{\bf Theorem 1.} {\it Under the hypothesis {\rm (H3)},

a) the limit
\begin{equation}\widetilde{\cal E}[\varphi,\chi]=\lim_n  \frac{\alpha_n}{2}\mathbb{E}[(\varphi(Y_n)-\varphi(Y))(\chi(Y_n)-\chi(Y)]\qquad \varphi, \chi\in{\cal D}\end{equation}
defines a closable positive bilinear form whose smallest closed extension is denoted $({\cal E},\mathbb{D})$.

b) $({\cal E},\mathbb{D})$ is a Dirichlet form

c) $({\cal E},\mathbb{D})$ admits a square field operator $\Gamma$ satisfying $\forall \varphi,\chi\in{\cal D}$
\begin{equation}
\Gamma[\varphi]=\widetilde{A}[\varphi^2]-2\varphi\widetilde{A}[\varphi]
\end{equation}
\begin{equation}\mathbb{E}_Y[\Gamma[\varphi]\chi]=\lim_n\alpha_n\mathbb{E}[(\varphi(Y_n)-\varphi(Y))^2(\chi(Y_n)+\chi(Y))/2]\end{equation}
\indent d) $({\cal E},\mathbb{D})$ is local if and only if $\forall \varphi\in{\cal D}$
\begin{equation}\lim_n \alpha_n\mathbb{E}[(\varphi(Y_n)-\varphi(Y))^4]=0\end{equation}
this condition is equivalent to
$\quad\exists\lambda> 2\quad\lim_n\alpha_n\mathbb{E}[|\varphi(Y_n)-\varphi(Y)|^\lambda]=0.$

e) If the form $({\cal E},\mathbb{D})$  is local, then the} principle of asymptotic error calculus
{\it is valid on 
$\widetilde{\cal D}=\{F(f_1,\ldots,f_p)\;:\;f_i\in{\cal D},\;\;F\in{\cal C}^1(\mathbb{R}^p,\mathbb{R})\}$
i.e.\\

$
\lim_n\alpha_n\mathbb{E}[(F(f_1(Y_n),\ldots,f_p(Y_n))-F(f_1(Y),\ldots,f_p(Y))^2]$\hfill

\hfill$=\mathbb{E}_Y[\sum_{i,j=1}^p F^\prime_i(f_1,\ldots,f_p)F^\prime_j(f_1,\ldots,f_p)\Gamma[f_i,f_j]].$}\\

An operator $B$ from ${\cal D}$ into $L^2(\mathbb{P}_Y)$ will be said to be a {\it first order operator} if it satisfies
$$B[\varphi\chi]=B[\varphi]\chi+\varphi B[\chi]\qquad\forall\varphi,\chi\in{\cal D}$$
\vspace{-.1cm}

\noindent{\bf Theorem 2.} {\it Under} (H1) {\it to} (H4). {\it If there is a real number $p\geq 1$ s.t.
$$\lim_n\alpha_n\mathbb{E}[(\varphi(Y_n)-\varphi(Y))^2|\psi(Y_n)-\psi(Y)|^p]=0\quad\forall\varphi,\psi\in{\cal D}$$
then $\Abar$ is first order.}

In particular, if the Dirichlet form is local, by the d) of theorem 1, the operator $\Abar$ is first order.\\

\noindent{\large\textsf{I.2. }} {\bf  Girsanov-type theorem for Dirichlet forms.}

An error structure is a probability space $(\Omega,\mathcal{A}, \mathbb{P})$ equipped with a local Dirichlet form with domain $\mathbb{D}$ dense in $L^2(\Omega,\mathcal{A}, \mathbb{P})$ admitting a square field operator $\Gamma$, see Bouleau [2]. We denote $\mathcal{D}A$ the domain of the associated generator.\\

\noindent{\bf Theorem 3.} {\it Let $(\Omega,\mathcal{A}, \mathbb{P}, \mathbb{D},\Gamma)$ be an error structure. Let be $f\in\mathbb{D}\cap L^\infty$ such that $f> 0$, $\mathbb{E}f=1$. We put $\mathbb{P}_1=f.\mathbb{P}$.

a) The bilinear form $\mathcal{E}_1$ defined on $\mathcal{D}A\cap L^\infty$ by
\begin{equation}\mathcal{E}_1[u,v]=-\mathbb{E}\left[fvA[u]+\frac{1}{2}v\Gamma[u,f]\right]\end{equation}
is closable in $L^2(\mathbb{P}_1)$ and satisfies for $u,v\in\mathcal{D}A\cap L^\infty$
\begin{equation}\mathcal{E}_1[u,v]=-\langle A_1u,v\rangle =-\langle u,A_1v\rangle =\frac{1}{2}\mathbb{E}[f\Gamma[u,v]]\end{equation}
where $A_1[u]=A[u]+\frac{1}{2f}\Gamma[u,f]$.

b) Let $(\mathbb{D}_1,\mathcal{E}_1)$ be the smallest closed extension of $(\mathcal{D}A\cap L^\infty,\mathcal{E}_1)$. Then $\mathbb{D}\subset\mathbb{D}_1$, $\mathcal{E}_1$ is local and admits a square field operator $\Gamma_1$, and 
$$\Gamma_1=\Gamma\quad\mbox{\rm on}\quad \mathbb{D}$$
in addition $\mathcal{D}A\subset\mathcal{D}A_1$ and 
$A_1[u]=A[u]+\frac{1}{2f}\Gamma[u,f]$  for all $u\in\mathcal{D}A$.}\\

\noindent{\bf Proof.} 1) First, using that the resolvent operators are bounded operators sending $L^\infty$ into $\mathcal{D}A\cap L^\infty$, we see that $\mathcal{D}A\cap L^\infty$ is dense in $\mathbb{D}$ (equipped with the usual norm $(\|.\|_{L^2}^2+\mathcal{E}[.])^{1/2}$), hence also dense in $L^2(\mathbb{P}_1)$.

2) Using that $\mathbb{D}\cap L^\infty$ is an algebra, for $u,v\in\mathcal{D}A\cap L^\infty$ we have
$$\mathcal{E}_1[u,v]=-\mathbb{E}[fvA[u]+\frac{1}{2}v\Gamma[u,f]]=\frac{1}{2}\mathbb{E} [\Gamma[fv,u]-v\Gamma[u,f]]=\frac{1}{2}\mathbb{E}[f\Gamma[u,v]].$$
So, defining $A_1$ as in the statement, we have $\forall u,v\in\mathcal{D}A\cap L^\infty$
$$\mathcal{E}_1[u,v]=-\mathbb{E}_1[vA_1u]=-\mathbb{E}_1[uA_1v].$$
The operator $A_1$ is therefore symmetric on $\mathcal{D}A\cap L^\infty$ under $\mathbb{P}_1$. Hence the form $\mathcal{E}_1$ defined on $\mathcal{D}A\cap L^\infty$ is closable, (Fukushima {\it et al.} [5], condition 1.1.3 p 4).

3) Let $(\mathbb{D}_1,\mathcal{E}_1)$ be the smallest closed extension of $(\mathcal{D}A\cap L^\infty,\mathcal{E}_1)$. Let be $u\in\mathbb{D}$ and $u_n\in\mathcal{D}A\cap L^\infty$, with $u_n\rightarrow u$ in $\mathbb{D}$. Using $\mathcal{E}_1[u_n-u_m]\leq \|f\|_\infty\mathcal{E}[u_n-u_m]$ and the closedness of $\mathcal{E}_1$ we get $u_n\rightarrow u$ in $\mathbb{D}_1$, hence $\mathbb{D}\subset\mathbb{D}_1$. Now by usual inequalities we see that $\Gamma[u_n]$ is a Cauchy sequence in $L^1(\mathbb{P}_1)$ and that the limit $\Gamma_1[u]$ does not depend on the particular sequence $(u_n)$ satisfying the above condition. Then following Bouleau [2], Chap. III \S2.5 p.38, the functional calculus extends to $\mathbb{D}_1$, the axioms of error structures are fulfilled for $(\Omega,\mathcal{A}, \mathbb{P}_1, \mathbb{D}_1,\Gamma_1)$ and this gives with usual arguments the b) of the statement.\hfill$\Box$\\

\noindent{\large\textsf{I.3. }} {\bf Rajchman measures.}

In the whole paper, if $x$ is a real number, $[x]$ denotes the entire part of $x$ and $\{x\}=x-[x]$ the fractional part.\\

\noindent{\bf Definition 1.} {\it A measure $\mu$ on the torus $\mathbb{T}^1$ is said to be Rajchman if 
$$\hat{\mu}=\int_{\mathbb{T}^1}e^{2i\pi nx}\,d\mu(x)\rightarrow 0\quad\quad\mbox{when }|n|\uparrow\infty.$$}

The set of Rajchman measures $\mathcal{R}$ is a band : if $\mu\in\mathcal{R}$ and if $\nu\ll|\mu|$ then $\nu\in\mathcal{R}$, cf.  Rajchman [18] [19], Lyons [15].\\

\noindent{\bf Lemma.} {\it Let $X$ be a real random variable and let $\Psi_X(u)=\mathbb{E}e^{iuX}$ be its characteristic function. Then 
$$\lim_{|u|\rightarrow\infty}\Psi_X(u)=0\quad\Longleftrightarrow\quad\mathbb{P}_{\{X\}}\in\mathcal{R}.$$}

\noindent{\bf Proof.} a) If $\lim_{|u|\rightarrow\infty}\Psi_X(u)=0$ then  $\Psi_X(2\pi n)=
(\mathbb{P}_{\{X\}})\hat{}\,(n)\rightarrow 0.$

b) Let $\rho$ be a probability measure on $\mathbb{T}^1$ s.t. $\rho\in\mathcal{R}$. From
$$e^{2i\pi ux}=e^{2i\pi[u]x}\sum_{p=0}^\infty \frac{((u-[u])2i\pi x)^p}{p!}$$ we have
$$\int e^{2i\pi ux}\rho(dx)=\sum_{p=0}^\infty \frac{((u-[u])2i\pi )^p}{p!}a_p([u])$$with
$a_p(n)=\int x^pe^{2i\pi nx}\rho(dx)$ hence
$|a_p(n)|\leq 1$ and $\lim_{|n|\rightarrow\infty}a_p(n)=0$ since $x^p\rho(dx)\in\mathcal{R}$, so 
$$\lim_{|u|\rightarrow\infty}\int e^{2i\pi ux}\rho(dx)=0.$$
Now if $\mathbb{P}_{\{X\}}\in\mathcal{R}$, since $1_{\{x\in[p,p+1[\}}.\mathbb{P}_{\{X\}}\ll\mathbb{P}_{\{X\}}$ we have
$$\lim_{|u|\rightarrow\infty}\mathbb{E}[e^{2i\pi uX}]=\lim_{|u|\rightarrow\infty}\sum_p \mathbb{E}[e^{2i\pi uX}1_{\{X\in[p,p+1[\}}]$$ which goes to zero by dominated convergence.\hfill$\Box$\\

A probability measure on $\mathbb{R}$ satisfying the conditions of the lemma  will be called Rajchman. \\

\noindent{\bf Examples.} Thanks to the Riemann-Lebesgue lemma, absolutely continuous measures are in $\mathcal{R}$. 
It follows from the lemma that if a measure $\nu$ satisfies $\nu\star\cdots\star\nu\in\mathcal{R}$ then $\nu\in\mathcal{R}$.
There are singular Rajchman measures,  cf. Kahane and Salem [13].\\

The preceding definitions and properties extend to $\mathbb{T}^d$ : a measure $\mu$ on $\mathbb{T}^d$ is said to be in $\mathcal{R}$ if $\hat{\mu}(k)\rightarrow0$ as $k\rightarrow\infty$ in $\mathbb{Z}^d$. The set of measures in $\mathcal{R}$ is a band. If $X $ is $\mathbb{R}^d$-valued, $\lim_{|u|\rightarrow\infty}\mathbb{E}e^{i\langle u,X\rangle }=0$ is equivalent to $\mathbb{P}_{\{X\}}\in\mathcal{R}$ where $\{x\}=(\{x_1\},\ldots,\{x_d\})$.\\

\noindent{\Large\textsf{II. Finite dimensional cases.}}\\

In the whole article $\stackrel{d}{\Longrightarrow}$ denotes the convergence in law, i.e. the convergence of the probability laws on bounded continuous functions. The {\it arbitrary functions principle} may be stated as follows:\\

\noindent{\bf Proposition 1.} {\it Let $X,Y,Z$ be random variables with values in $\mathbb{R}$, $\mathbb{R}$, and $\mathbb{R}^m$ resp. Then 
\begin{equation}(\{nX+Y\},X,Y,Z)\quad\stackrel{d}{\Longrightarrow}\quad(U,X,Y,Z)\end{equation}
where $U$ is uniform on the unit interval independent of $(X,Y,Z)$, if and only if $\mathbb{P}_X$ is Rajchman.}\\

\noindent{\bf Proof.} If $\mu$ is a probability measure on $\mathbb{T}^1\times\mathbb{R}^m$, let us put 
$$\hat{\mu}(k,\zeta)=\int e^{2i\pi kx+\langle \zeta,y\rangle }\mu(dx,dy),$$ then $\mu_n\stackrel{d}{\Longrightarrow}\mu$ iff $\hat{\mu}_n(k,\zeta)\rightarrow\hat{\mu}(k,\zeta)$ $\forall k\in\mathbb{Z}$, $\forall\zeta\in\mathbb{R}^m$.

a) If $\mathbb{P}_X\in\mathcal{R}$ 
$$\hat{\mathbb{P}}_{(\{nX+Y\},X,Y,Z)}(k,\zeta_1,\zeta_2,\zeta_3)=\mathbb{E}[\exp\{2i\pi k(nX+Y)+i\zeta_1X+i\zeta_2Y+i\langle \zeta_3,Z\rangle \}]$$
$$=\int e^{2i\pi knx}f(x)\mathbb{P}_{\{X\}}(dx)$$
with $f(x)=\mathbb{E}[\exp\{2i\pi kY+i\zeta_1X+i\zeta_2Y+i\langle \zeta_3,Z\rangle \}|\{X\}=x]$. The fact that $f.\mathbb{P}_{\{X\}}\in\mathcal{R}$ gives the result.

b) Conversely, taking $(k,\zeta_1,\zeta_2,\zeta_3)= (1,0,-2\pi, 0)$ gives  $\hat{\mathbb{P}}_{\{X\}}(n)\rightarrow 0$ i.e. $\mathbb{P}_X\in\mathcal{R}$.\hfill$\Box$\\

 Let us  suppose now that $Y$ is an $\mathbb{R}^d$-valued random variable, measured with an equidistant graduation corresponding to an orthonormal rectilinear coordinate system, and estimated to the nearest graduation component by component. Thus we put
$$Y_n=Y+\frac{1}{n}\theta(nY)$$ with $\theta(y)=(\frac{1}{2}-\{y_1\},\cdots,\frac{1}{2}-\{y_d\})$. Let us emphasize that $Y_n$ is a deterministic function of $Y$.\\

\noindent{\bf Theorem 4.} {\it a) If $\mathbb{P}_Y$ is Rajchman and if $X$ is $\mathbb{R}^m$-valued

\begin{equation}(X,n(Y_n-Y))\quad\stackrel{d}{\Longrightarrow}\quad(X,(V_1,\ldots,V_d))\end{equation}
where the $V_i$'s are independent identically distributed uniformly distributed on $(-\frac{1}{2},\frac{1}{2})$ and independent of $X$.

\noindent For all $\varphi\in\mathcal{C}^1\cap lip(\mathbb{R}^d)$
\begin{equation}
(X,n(\varphi(Y_n)-\varphi(Y)))\quad\stackrel{d}{\Longrightarrow}\quad(X,\sum_{i=1}^dV_i\varphi^\prime_i(Y))
\end{equation}
\begin{equation}
n^2\mathbb{E}[(\varphi(Y_n)-\varphi(Y))^2Ü|Y\!=\!y]\rightarrow\frac{1}{12}\sum_{i=1}^d\varphi^{\prime 2}_i(y)\qquad\mbox{ in }L^1(\mathbb{P}_Y)
\end{equation}
in particular

\begin{equation}
n^2\mathbb{E}[(\varphi(Y_n)-\varphi(Y))^2]\rightarrow\mathbb{E}_Y[\frac{1}{12}\sum_{i=1}^d\varphi^{\prime 2}_i(y)].
\end{equation}
\indent b) If $\varphi$ is of class $\mathcal{C}^2$, the conditional expectation $n^2\mathbb{E}[\varphi(Y_n)-\varphi(Y)|Y=y]$ possesses a version $n^2(\varphi(y+\frac{1}{n}\theta(ny))-\varphi(y))$ independent of the probability measure $\mathbb{P}$ which converges in the sense of distributions to the function $\frac{1}{24}\bigtriangleup\varphi$.

c) If $\mathbb{P}_Y\ll dy $ on $\mathbb{R}^d$, $\forall\psi\in L^1([0,1])$
\begin{equation}(X,\psi(n(Y_n-Y)))\quad\stackrel{d}{\Longrightarrow}(X,\psi(V)).\end{equation}

d) We consider the bias operators on the algebra $\mathcal{C}^2_b$ of bounded functions with bounded derivatives up to order 2 with the sequence $\alpha_n=n^2$. If $\mathbb{P}_Y\in\mathcal{R}$ and if one of the following condition is fulfilled

i) $\forall i=1,\ldots,d$ the partial derivative $\partial_i\mathbb{P}_Y$ in the sense of distributions is a measure $\ll\mathbb{P}_Y$ of the form $\rho_i\mathbb{P}_Y$ with $\rho_i\in L^2(\mathbb{P}_Y)$,

ii) $\mathbb{P}_Y=h1_{G}\frac{dy}{|G|}$  with $G$ open set, $h\in H^1\cap L^\infty(G)$, $h> 0$,

\noindent then hypotheses {\rm (H1)} to {\rm (H4)} are satisfied and 
$$\begin{array}{rl}
\overline{A}[\varphi]&=\frac{1}{24}\bigtriangleup\varphi\\
\widetilde{A}[\varphi]&=\frac{1}{24}\bigtriangleup\varphi+\frac{1}{24}\sum\varphi^\prime_i\rho_i\qquad\mbox {case i)}\\
\widetilde{A}[\varphi]&=\frac{1}{24}\bigtriangleup\varphi+\frac{1}{24}\frac{1}{h}\sum h^\prime_i\varphi^\prime_i\qquad\mbox {case ii)}\\
\Gamma[\varphi]&=\frac{1}{12}\sum \varphi^{\prime 2}_i.
\end{array}
$$}

\noindent{\bf Proof.}  The argument for relation (8) is similar to the one dimensional case stated in proposition 1. The relation (9) comes from the Taylor expansion 
$\varphi(Y_n)-\varphi(Y)=$

$=\sum_{i=1}^d (Y_{n,i}-Y_i)\int_0^1\varphi^\prime_i(Y_{n,1},\ldots,Y_{n,i-1},Y_i+t(Y_{n,i}-Y_i),Y_{i+1},\ldots,Y_d)\,dt$

\noindent and the convergence
$$(X,\sum_i\theta(nY_i)\varphi^\prime_i(Y))\quad\stackrel{d}{\Longrightarrow}\quad(X,\sum_i\varphi^\prime_i(Y)V_i)$$
thanks to (8) and the following approximation in $L^1$ 
$$\mathbb{E}\left|\sum_i\theta(nY_i)\varphi^\prime_i(Y)-\sum_i\theta(nY_i)\int_0^1\varphi^\prime_i(\ldots,Y_i+t(Y_{n,i}-Y_i),\ldots)dt\right|\rightarrow0.$$
To prove the formulas (10) and (11) let us remark that

\begin{flushleft}
$n^2\mathbb{E}[(\varphi(Y_n)-\varphi(Y)^2|Y=y]=$
\end{flushleft}
$$=\mathbb{E}\left[\left|\sum_i\theta(nY_i)\int_0^1\varphi^\prime_i(\ldots,Y_i+t(Y_{n,i}-Y_i),\ldots)dt\right|^2|Y=y\right]$$
$$=\left|\sum_i\theta(ny_i)\int_0^1\varphi^\prime_i(y_1+\frac{1}{n}\theta(ny_1),\ldots,y_i+t\frac{1}{n}\theta(ny_i),\ldots)dt\right|^2\quad\mathbb{P}_Y{\mbox{-a.s.}}$$
each term $(\theta(ny_i)\int_0^1\varphi^\prime_i(\ldots)dt)^2$ converges to $\int\theta^2\varphi^{\prime 2}_i(y)=\frac{1}{12}\varphi^{\prime 2}_i$ in $L^1$ and each term $\theta(ny_i)\theta(ny_j)\int_0^1\ldots\int_0^1\ldots$ goes to zero in $L^1$ what proves the part a) of the statement.

The part b) is obtained following the same lines with a Taylor expansion up to second order and an integration by part thanks to the fact that $\varphi$ is now supposed to be $\mathcal{C}^2$.

In order to prove c) let us suppose first that $\mathbb{P}_Y=1_{[0,1]^d}.dy$. Considering a sequence of functions $\psi_k\in\mathcal{C}_b$ tending to $\psi$ in $L^1$ we have the bound 
$$
\begin{array}{l}
|\mathbb{E}[e^{i\langle u,X\rangle }e^{iv\psi(\theta(nY))}]-\mathbb{E}[e^{i\langle u,X\rangle }e^{iv\psi_k(\theta(nY))}]|\\
\leq |v|\int|\psi(\theta(ny))-\psi_k(\theta(ny))|dy\\
=|v|\sum_{p_1=0}^{n-1}\cdots\int_{p_1}^{p_1+1}\cdots|\psi(\theta(ny_1)\ldots)-\psi_k(\theta(ny_1)\ldots)|dy_1\ldots dy_d\\
=|v|\sum\cdots\sum\int\cdots\int|\psi(\theta(x_1),\ldots)-
\psi_k(\theta(x_1),\ldots)|\frac{dx_1}{n}\cdots\frac{dx_d}{n}\\
=|v|\|\psi-\psi_k\|_{L^1}.
\end{array}
$$
This yields (12) in this case. Now if $\mathbb{P}_Y\ll dy$ then $\mathbb{P}_{\{Y\}}\ll dy $ on $[0,1]^d$ and the weak convergence under $dy$ on $[0,1]^d$ implies the weak convergence under $\mathbb{P}_{\{Y\}}$ what yields the result.

In d) the point i) is proved by the approach already used in Bouleau [3] consisting of proving  that hypothesis (H3) is fulfilled by displaying the operator $\widetilde{A}$ thanks to an integration by parts. The point ii) is an application of the Girsanov-type theorem 3.\hfill$\Box$\\

\noindent{\bf Remarks.} 1) About the relations (9) (10) (11), let us note that with respect to the form $$\mathcal{E}[\varphi]=\frac{1}{24}\mathbb{E}_Y \sum_i\varphi^{\prime 2}_i$$ when it is closable, the random variable $\sum_iV_i\varphi^\prime_i$ appears to be {\it a gradient} :
if we put $\varphi^\#=\sum_iV_i\varphi^\prime_i$ then a we have
$$\mathbb{E}[\varphi^{\# 2}]=\frac{1}{12}\sum_i\varphi^{\prime 2}_i=\Gamma[\varphi]$$
the square field operator associated to $\mathcal{E}$. We will find this phenomenon again on the Wiener space.

2) {\it Approximation to the nearest graduation, by excess, or by default.}
 When the approximation is done to the nearest graduation, on the algebra $\mathcal{C}^2_b$ the four bias operators are obtained in theorem 4 with the sequence $\alpha_n=n^2$, (with $\alpha_n=n$ the four bias operators would be zero).

We would obtain a quite different result with an approximation by default or by excess because of the dominating effect of the shift.

If the random variable $Y$ is approximated by default by $Y_n^{(d)}=\frac{[nY]}{n}$ then 
$$n(Y_n^{(d)}-Y)\stackrel{d}{\Longrightarrow} -U\quad\mbox{and}\quad\mathbb{E}[n(Y_n^{(d)}-Y)]\rightarrow -\frac{1}{2}$$
as soon as $Y$ is say bounded. With this approximation, if we do not erase the shift down proportional to $-\frac{1}{2n}$, and if we take $\alpha_n=n$ we obtain  first order bias operators without diffusion :
$\overline{A}[\varphi]=-\frac{1}{2}\varphi^\prime=-\underline{A}[\varphi]$ and $\widetilde{A}=0$. The same happens of course with the approximation by excess.\\

3) {\it Extension to more general graduations.}
Let $Y$ be an $\mathbb{R}^d$-valued random variable approximated by $Y_n=Y+\xi_n(Y)$ with a sequence $\alpha_n\uparrow\infty$ on the algebra $\mathcal{D}=\mathcal{L}\{e^{\langle u,x\rangle },\;u\in\mathbb{R}^d\}$, the function $\xi_n$ satisfying 
$$
(\ast)\left\{\begin{array}{l}
\alpha_n\mathbb{E}[|\xi_n|^3(Y)]\rightarrow 0\\
\\
\alpha_n\mathbb{E}[\varphi(Y)\langle u,\xi_n(Y)\rangle ^2]\rightarrow\mathbb{E}_Y[\varphi . u^\ast\underline{\underline{\gamma}}u]\qquad\forall\varphi\in\mathcal{D},\forall u\in\mathbb{R}^d\\
\mbox{ with }\gamma_{ij}\in L^\infty(\mathbb{P}_Y)\mbox{ and }\frac{\partial \gamma_{ij}}{\partial x_j}\mbox{ in distributions sense }\in L^2(\mathbb{P}_Y)\\
\\
\alpha_n\mathbb{E}[\varphi(Y)\langle u,\xi_n(Y)\rangle ]\rightarrow 0\quad \forall\varphi\in\mathcal{D}.
\end{array}
\right.
$$
Under these hypotheses we have\\

\noindent{\bf Theorem 4bis.} {\it a) {\rm (H1)} is satisfied and 
$\overline{A}[\varphi]=\frac{1}{2}\sum_{ij}\gamma_{ij}\frac{\partial^2\varphi}{\partial x_i\partial x_j}.$

b) If for $i=1,\ldots,d$, the partial derivative $\partial_i\mathbb{P}_Y$ in the sense of distributions is a bounded measure of the form $\rho_i\mathbb{P}_Y$ with $\rho_i\in L^2(\mathbb{P}_Y)$ then assumptions {\rm (H1)} to {\rm (H4)} are fulfilled and $\forall\varphi\in\mathcal{D}$
$$\widetilde{A}[\varphi]=\frac{1}{2}\sum_{ij}\gamma_{ij}\frac{\partial^2\varphi}{\partial x_i\partial x_j}+\sum_i(\sum_j(\frac{\partial \gamma_{ij}}{\partial x_j}+\gamma_{ij}\rho_j))\frac{\partial \varphi}{\partial x_i}$$
the square field operator is 
$\Gamma[\varphi]=\sum_{ij}\gamma_{ij}\frac{\partial \varphi}{\partial x_i}\frac{\partial \varphi}{\partial x_j}.$}\\

\noindent{\bf Proof.} The argument is simple thanks to the choice of  the algebra $\mathcal{D}$ and consists of elementary Taylor expansions  to prove the existence of the bias operators. Then theorem 1 applies.\hfill$\Box$\\

\noindent{\bf Historical comment.}

In his intuitive version, the idea underlying the arbitrary functions method is ancient. The historian J. von Plato [16] dates it back to a book of J. von Kries [12]. We find indeed in this philosophical treatise the idea that if a roulette had equal and infinitely small black and white cases, then there would be an equal probability to fall on a case or on the neighbour one, hence by addition an equal probability to fall either on black or on white. But no precise proof was given. The idea remains at the common sense level. 

A mathematical argument for the fairness of the roulette and for the equi-distribu\-tion of other mechanical systems (little planets on the Zodiac) was proposed by H. Poincar\'e  in his course on probability published in 1912 ([17], Chap. VIII \S 92 and especially \S 93). In present language, Poincar\'e shows the weak convergence of $tX+Y \mbox{mod }  2\pi$ when $t\uparrow \infty$ to the uniform law on $(0,2\pi)$ when the pair $(X,Y)$ has a density. He uses the characteristic functions. His proof supposes the density  be $\mathcal{C}^1$ with bounded derivative in order to perform an integration by parts, but the proof would extend to the general absolutely case if we were using instead the Riemann-Lebesgue lemma. 

The question is then developed  without major changes by several authors, E. Borel [1] (case of continuous density), M. Fr\'echet [4] (case of Riemann-integrable density), B. Hostinski [9] [10] (bidimensional case) and is tackled anew by E. Hopf [6] , [7] and [8] with the more general point of view of asymptotic behaviour of dissipative dynamical systems. Hopf has shown that these phenomena are related to mixing and belong to the framework of ergodic theory. \\

\noindent{\Large\textsf{III. Infinite dimensional extensions of the arbitrary functions principle.}}\\

\noindent{\bf III.1. Rajchman type martingales.}

Let $(\mathcal{F}_t)$ be a right continuous filtration on $(\Omega, \mathcal{A}, \mathbb{P})$ and $M$ be a continuous local $(\mathcal{F}_t,\mathbb{P})$-martingale nought at zero. $M$ will be said to be Rajchman if the measure $d\langle M,M\rangle_s$ restricted to compact intervals belongs to $\mathcal{R}$ almost surely. 
 We will show that the method followed by Rootz\'en [21] extends to Rajchman martingales and provides the following\\

\noindent{\bf Theorem 5.} {\it Let $M$ be a continuous local martingale which is Rajchman and s.t. $\langle M,M\rangle_\infty=\infty$.

Let $f$ be a bounded Riemann-integrable periodic function with unit period on $\mathbb{R}$ s.t. $\int_0^1f(s)ds=0$.
 Then for any random variable $X$
\begin{equation}(X,\int_0^.f(ns)\,dM_s)\quad\stackrel{d}{\Longrightarrow}\quad(X,W_{\|f\|^2\langle M,M\rangle_.}),\end{equation}
the weak convergence is understood on $\mathbb{R}\times\mathcal{C}([0,1])$ and  $W$ is an independent standard Brownian motion.}\\

Before proving the theorem, let us remark that it shows that the random measure $dM_s$ behaves in some sense like a Rajchman measure. Indeed if $\mathbb{P}_Y\in\mathcal{R}$ we have 
$$\int_{-\infty}^yg(nx)\mathbb{P}_Y(dx)\rightarrow \int_0^1g(x)dx\int_{-\infty}^y\mathbb{P}_Y(dx)$$
as soon as $g$ is periodic with unit period, Riemann-integrable and bounded. Now applying the theorem to the Brownian motion gives the similar relation
$$\int_0^t f(ns)\,dB_s\quad\stackrel{d}{\Longrightarrow}\quad(\int_0^1f^2(s)ds)^{1/2}\int_0^tdW_s.$$
{\bf Proof.} We consider the local martingale $N_t=\int_0^tf(ns)dM_s$.

a) In order to be sure that $\langle N,N\rangle_\infty=\infty$, we change $N_t$ into $\tilde{N_t}=\int_0^tf_n(s)dM_s$ with $f_n(s)=f(ns)$ for $s\in[0,1)$, $f_n(s)=0$ for $s\in[1,n]$ and $f_n(s)=1$ for $t> n$. We put $S_n(t)=\inf\{s\,:\langle\tilde{N},\tilde{N}\rangle_s\;> \;t\}$.

b) We want to show
\begin{equation}\qquad\mathbb{E}[\xi F(\tilde{N}_{S_n})]\rightarrow\mathbb{E}[\xi F(W)]\quad\forall\xi\in L^1(\mathbb{P})\quad\forall F\in \mathcal{C}_b([0,1]).\end{equation}
It is enough to consider the case $\xi> 0$, $\mathbb{E}\xi=1$, and $\xi$ may be supposed to be $\mathcal{F}_T$-measurable for a deterministic time $T$ large enough. Let be $\tilde{\mathbb{P}}=\xi.\mathbb{P}$ and $D(t)=\mathbb{E}[\xi|\mathcal{F}_t]$. The process
$$\tilde{M_t}=M_t-\int_0^t D^{-1}(s)d\langle M,D^c\rangle_s$$ is a continuous local martingale under $\tilde{\mathbb{P}}$.  Therefore $\int_0^{S_n(t)}f_n(s)\,d\tilde{M}_s$ is a Brownian motion under $\tilde{\mathbb{P}}$ (Revuz and Yor [20]  p.313 theorem 1.4 and p 173). Writing
$$\int_0^{S_n(t)}f_n(s)dM_s=\int_0^{S_n(t)}f_n(s)d\tilde{M}_s+\int_0^{S_n(t)}\frac{f_n(s)}{D(s)}d\langle M,D^c\rangle_s$$ and noting that $d\langle M,D^c\rangle_s$ vanishes on $]T,\infty[$, in order to show $(14)$ it suffices to show
$$\sup_{0\leq t\leq T}\left|\int_0^t\frac{f_n(s)}{D(s)}\,d\langle M,D^c\rangle_s\right|\rightarrow 0\quad \mbox{ a.s. when } n\rightarrow\infty$$
hence to show
$$\sup_{0\leq t\leq 1}\left|\int_0^t\frac{f(ns)}{D(s)}\,d\langle M,D^c\rangle_s\right|\rightarrow 0\quad \mbox{ a.s. when } n\rightarrow\infty$$
and, because $M$ is Rajchman this comes from the following lemma :\\

\noindent{\bf Lemma.} {\it Let $f$ be as in the statement of the theorem, then $\forall\mu\in\mathcal{R}$
$$\sup_{0\leq t\leq 1}\left|\int_0^tf(ns)\mu(ds)\right|\rightarrow 0\quad \mbox{ as } n\rightarrow\infty.$$}

\noindent{\bf Proof.} We have 
$$\int_0^t f(ns)\mu(ds)\rightarrow\int_0^1f(s)ds\int_0^t\mu(ds)=0.$$
Since $f$ is bounded, the functions $\int_0^tf(ns)\mu(ds)$ are equi-continuous and the result follows from Ascoli theorem.\hfill$\Box$\\

c) This proves the following stable convergence
$$(X,\int_0^{T_n(.)}f(ns)dM_s)\quad\stackrel{d}{\Longrightarrow}\quad(X,W_.)$$ and by the fact that the limit 
$$\int_0^t f^2(ns)d\langle M,M\rangle_s\rightarrow\int_0^1f^2(s)ds\langle M,M\rangle_t$$ is a continuous process, this gives the announced result.\hfill$\Box$\\

\noindent{\bf Remark.} If $\int_0^1f(s)ds\neq 0$, then keeping the other hypotheses unchanged, we obtain
$$(X,\int_0^.f(ns)dM_s)\quad\stackrel{d}{\Longrightarrow}\quad\left(X,(\int_0^1f(s)ds)M_.+(\int_0^1(f-\!\int_0^1f)^2)^{1/2}W_{\langle M,M\rangle_.}\right).$$
\vspace{.2cm}

We study now the induced limit quadratic form when the martingale $M$ is approximated by the martingale $M^n_t=M_t+\int_0^t\frac{1}{n}f(ns)dM_s$. The notation is the same as in the preceding section and $f$ satisfies the same hypotheses as in theorem 5.\\

\noindent{\bf Theorem 6.} {\it Let $M$ be a Rajchman martingale s.t. $M_1\in L^2$ and $\eta$, $\zeta$ bounded adapted processes. Then
$$
\begin{array}{c}
\displaystyle n^2\mathbb{E}\left[(\exp\{i\int_0^1\eta_sdM^n_s\}-\exp\{i\int_0^1\eta_sdM_s\})(\exp\{i\int_0^1\zeta_sdM^n_s\}-\exp\{i\int_0^1\zeta_sdM_s\})\right]\\
\\
\displaystyle\rightarrow-\mathbb{E}\left[\exp\{i\int_0^1(\eta_s+\zeta_s)dM_s\}\int_0^1\eta_s\zeta_s\,d\langle M,M\rangle_s\right]\int_0^1f^2(s)ds.
\end{array}
$$}

\noindent{\bf Proof.} By the fundamental formula of calculus (finite increments formula), the first term in the statement may be written
$$-\mathbb{E}[\exp\{i\int_0^1(\eta_s+\zeta_s)dM_s\}\int_0^1\eta_s f(ns)dM_s\int_0^1\zeta_sf(ns)dM_s]+o(1)$$ therefore, thanks to theorem 5, the statement is a consequence of the following lemma :\\

\noindent{\bf Lemma.} {\it Suppose $\mathbb{E}M_1^2<\infty$ and $\eta$ adapted and bounded, then the random variables $\int_0^1\eta_sf(ns)dM_s$ are uniformly integrable.}\\

\noindent{\bf Proof.} It suffices to remark that their $L^2$-norm is equal to $\mathbb{E}\int_0^1\eta_s^2f^2(ns)\,d\langle M,M\rangle_s$ hence uniformly bounded.\hfill$\Box$\\

\noindent{\bf III.2. Sufficient closability conditions on the Wiener space.}

The closability problem of the limit quadratic forms obtained in the preceding section, may be tackled with the tools available on the Wiener space.

Let us approximate the Brownian motion $(B_t)_{t\in[0,1]}$ by the process $B^n_t=B_t+\int_0^t\frac{1}{n}f(ns)\,dB_s$ where $f$ satisfies the same hypotheses as before.\\

\noindent{\bf Theorem 7.} {\it a) Let $\xi\in L^2([0,1])$, and let $X$ be a random variable defined on the Wiener space, i.e. a Wiener functional, then
\begin{equation}
\left(X,n(\exp\{i\int_0^1\xi dB^n\}-\exp\{i\int_0^1\xi dB\})\right)
\stackrel{d}{\Longrightarrow}\left(X,\|f\|_{L^2}(\exp\{i\int_0^1\xi dB\})^{\#}\right)
\end{equation}
here for any regular Wiener functional $Z$ we put $Z^\#(\omega,w)=\int_0^1D_sZ\,dW_s$, where $W$ is an independent Brownian motion.

b) 
\begin{equation}
n^2\mathbb{E}\left[(e^{i\xi.B^n}-e^{i\xi.B})^2\right]\rightarrow-\mathbb{E}[e^{2i\xi.B}]\int_0^1\xi^2ds\|f\|^2_{L^2}\end{equation}
on the algebra $\mathcal{L}\{e^{i\xi.B}\}$ the quadratic form $-\frac{1}{2}\mathbb{E}[e^{2i\xi.B}]\int_0^1\xi^2ds$ is closable, its closure is the Ornstein-Uhlenbeck form.}\\

\noindent{\bf Proof.} a) The first assertion comes easily from the similar result concerning Rajchman martingales using the fact that
$\int_0^1e^{i\alpha\int_0^1\frac{1}{n}f(ns)dB_s}d\alpha\rightarrow1$ in $L^p$ $p\in[1,\infty[$.

b) The obtained quadratic form is immediately recognized as the Ornstein-Uhlenbeck form which is closed. It follows that hypothesis (H3) is fulfilled and the symmetric bias operator is
$$\widetilde{A}[e^{i\int\xi dB}]=\left(-\frac{i}{2}\int\xi dB-\frac{1}{2}\int\xi^2 ds\right)e^{i\int \xi dB}.\hspace{2cm}\Box$$

If instead of the Wiener measure $m$, we consider the measure $m_1=h.m$ for an $h> 0$, $h\in\mathbb{D}_{ou}\cap L^\infty$ where $\mathbb{D}_{ou}$ $(=D^{2,1})$ denotes the domain of the Ornstein-Uhlenbeck form, we know by the Girsanov-type theorem 3 that the form 
$$-\frac{1}{2}\mathbb{E}_1[e^{2i\xi.B}\int_0^1\xi^2ds]$$ is closable, admits the same square field operator on $\mathbb{D}_{ou}$, and that its generator $A_1$ satisfies
$$A_1[\varphi]=\widetilde{A}[\varphi]+\frac{1}{2h}\Gamma_{ou}[\varphi,h]\quad\mbox{ for }\quad\varphi\in\mathcal{D}A_{ou}$$
Since the point a) of the theorem is still valid under $m_1$ because of the properties of stable convergence, the preceding theorem is still valid under $m_1$, the Dirichlet form being now
$$\mathcal{E}_1[\varphi]=\frac{1}{2}\mathbb{E}_1[\Gamma_{ou}[\varphi]]\quad\mbox{ for }\quad \varphi\in\mathbb{D}_{ou}.$$ 

\noindent{\bf Remark.} Let us come back to the general case of Rajchman martingales. If we suppose the Rajchman local martingale $M$ is in addition Gaussian, which is equivalent to suppose $\langle M,M\rangle$ deterministic, then  on the algebra $\mathcal{L}\{e^{i\int\xi dM}; \xi$  deterministic bounded $\}$ the limit quadratic form 
$$-\mathbb{E}[e^{i\int(\eta+\zeta)dM}\int_0^1\zeta_s\eta_sd\langle M,M\rangle_s]\|f\|^2_{L^2}$$
is closable, hence (H3) is satisfied.

Indeed, it suffices to exhibit the corresponding symmetric bias operator. But by the use of the calculus for Gaussian variables, it is easily seen that the operator defined by
$$\widetilde{A}[e^{i\int\xi dM}]=e^{i\int\xi dM}\left(-\frac{i}{2}\int\xi dM-\frac{1}{2}\int\xi^2\,d\langle M,M\rangle_s\right)\int f^2 ds$$ satisfies the required condition.$\hfill$$\Box$\\

\noindent{\bf III.3.  Approximation of the Ornstein-Uhlenbeck gradient.}

Let $m$ be the Wiener measure on $\mathcal{C}([0,1],\mathbb{R})$. Let $\theta$ be a real periodic function of period 1 such that  $\int_0^1\theta(s)ds=0$ et $\int_0^1\theta^2(s)ds=1$. We consider the transformation $R_n$ of the space $L^2_\mathbb{C}(m)$ defined by its action on the chaos:

\noindent if $X=\int_{s_1<\cdots<s_k}\hat{f}(s_1,\ldots,s_k)dB_{s_1}\ldots dB_{s_k}$ for $\hat{f}\in L^2_{sym}([0,1]^k,\mathbb{C})$,
$$R_n(X)=\int_{s_1<\cdots<s_k}\hat{f}(s_1,\ldots,s_k)e^{i\frac{1}{n}\theta(ns_1)}dB_{s_1}\ldots e^{i\frac{1}{n}\theta(ns_k)}dB_{s_k}.$$
Since $\|X\|^2_{L^2(m)}=\int_{s_1<\cdots<s_k}|\hat{f}|^2ds_1\ldots ds_k=\frac{1}{k!}\|\hat{f}\|^2_{L^2_{sym}}$,  $R_n$ is an isometry from $L^2_\mathbb{C}(m)$ into itself and $\forall \xi\in L^2_\mathbb{C}([0,1])$
$$R_n[e^{\int\xi dB_s-\frac{1}{2}\int\xi^2ds}]=e^{\int\xi e^{i\frac{1}{n}\theta(ns)}dB_s-\frac{1}{2}\int\xi^2e^{\frac{2i}{n}\theta(ns)}ds}$$
$$\|e^{\int\xi dB_s-\frac{1}{2}\int\xi^2ds}\|_{L^2_\mathbb{C}}=e^{\frac{1}{2}\int|\xi|^2 ds}.$$
From the relation
$$n(e^{\frac{i}{n}\sum_{p=1}^k\theta(ns_p)}-1)=i\sum_{p=1}^k\theta(ns_p)\int_0^1 e^{\alpha \frac{i}{n}\sum_p\theta(ns_p)}d\alpha$$
it follows that if $X$ belongs to $k$-th chaos
$$\|n(R_n(X)-X)\|^2_{L^2}\leq k^2\|X\|^2\|\theta\|^2_\infty$$
then, denoting $A$ the Ornstein-Uhlenbeck operator,  for $X\in \mathcal{D}(A)$
$$\|n(R_n(X)-X)\|_{L^2}\leq 2\|AX\|\|\theta\|_\infty$$ and we can state:\\

\noindent{\bf Theorem 8.} {\it If $X\in\mathcal{D}(A)$
$$(-in(R_n(X)-X),B)\quad\stackrel{d}{\Longrightarrow}\quad(X^\#,B)$$
with $X^\#=\int_0^1D_sX\,dW_s$ where $W$ is an independent Brownian motion.}\\

\noindent{\bf Proof.}  If $X$ belongs to the $k$-th chaos, expanding the exponential by its Taylor series gives 
$$n(R_n(X)-X) =i\int_{s_1<\cdots<s_k}\hat{f}(s_1,\ldots,s_k)\sum_{p=1}^k\theta(ns_p)dB_{s_1}\ldots dB_{s_k}+Q_n$$
with $\|Q_n\|^2\leq \frac{1}{4n}k^2\|\theta\|^2_\infty\|X\|^2$.

Now, since $\int_{s_1<\cdots<s_p<\cdots<s_k}h(s_1,\ldots,s_k)\theta(ns_p)dB_{s_1}\ldots dB_{s_p}\ldots dB_{s_k}$ converges stably to
$\int_{s_1<\cdots<s_p<\cdots<s_k}h(s_1,\ldots,s_k)dB_{s_1}\ldots dW_{s_p}\ldots dB_{s_k}$ we obtain that 
$$-in(R_n(X)-X)\stackrel{s}{\Longrightarrow}
\begin{array}{l}
\int_{t<s_2<\cdots<s_k}\hat{f}(t,s_2,\ldots,s_k)dW_tdB_{s_2}\ldots dB_{s_k}\\
+\int_{s_1<t<\cdots<s_k}\hat{f}(s_1,t,\ldots,s_k)dB_{s_1}dW_t\ldots dB_{s_k}\\
+\cdots\\
+\int_{s_1<\cdots<s_{k-1}<t}\hat{f}(s_1,\ldots,s_{k-1},t)dB_{s_1}\ldots dB_{s_{k-1}}dW_t
\end{array}
$$
which is equal to $\int D_s(X)dW_s=X^\#$.

For the general case, we approximate $X$ by $X_k$ for the norm $\mathbb{D}^{2,2}$ and  reasoning with the characteristic functions yields the result (see the proof of theorem 10 below). \hfill$\Box$\\

By the properties of the stable convergence, the convergence in law of theorem 8 still holds under  $\tilde{m}\ll m$.\\

\noindent{\bf Theorem 9.} {\it $\forall X\in\mathcal{D}(A)$
$$n^2\mathbb{E}[|R_n(X)-X|^2]\rightarrow2\mathcal{E}[X]$$
where $\mathcal{E}$ is the Dirichlet form associated with the Ornstein-Uhlenbeck operator.}\\

\noindent{\bf Proof.} As $R_n$ preserves the chaos and the expansion of $n(R_n(X)-X)$ on the chaos is dominated by that of $2\|\theta\|_\infty AX$, it suffices to argue when $X$ is in the $k$-th chaos.

Starting from 
$$n^2\mathbb{E}|R_n(X)-X|^2=n^2\int_{s_1<\cdots<s_k}|\hat{f}|^2\,|e^{\frac{i}{n}\sum_p\theta(ns_p)}-1|^2\,ds_1\ldots ds_k$$
expanding the exponential  and estimating the remainder we obtain $$\lim_n n^2\mathbb{E}|R_n(X)-X|^2=\frac{k}{k!}\int_{[0,1]^k}|\hat{f}|^2ds_1\ldots ds_k\int_0^1\theta^2 dt=k\|X\|^2$$
what gives the result.\hfill$\Box$\\

Following the same arguments, it is possible to show that the theoretical  and practical bias operators $\overline{A}$ and $\underline{A}$ defined on the algebra $\mathcal{L}\{e^{\int\xi dB}\;;\;\xi\in\mathcal{C}^1\}$ by
$$\begin{array}{c}
n^2\mathbb{E}[(R_n(X)-X)Y]=\langle \overline{A}X,Y\rangle _{L^2(m)}\\
n^2\mathbb{E}[(X-R_n(X))R_n(Y)]=\langle \underline{A}X,Y\rangle _{L^2(m)}
\end{array}
$$
exist and are equal to $A$.\\

\noindent{\bf III.4. Isometries on the Wiener space.}

Let us now consider a $d$-dimensional Brownian motion $(B_t)$. 

Let $t\mapsto M_t$ be a deterministic bounded measurable periodic map with period 1 with values in the space of $d\times d$ othogonal matrices such that $\int_0^1M_s ds=0$, (for instance a rotation of angle $2\pi t$). We denote still $m$ the Wiener measure law of $(B_t)$. The transformation $B_t\mapsto\int_0^tM_sdB_s$ induces an endomorphism $T_M$ isometric in $L^p(m), 1\leq p\leq \infty$. We put $M_n(s)=M_{ns}$ and $T_n=T_{M_n}$.\\

\noindent{\bf Theorem 10.} {\it Let $X$ be in $ L^1(m)$. Let be $\tilde{m}\ll m$, we have under $\tilde{m}$:
$$(T_n(X),B)\quad\stackrel{d}{\Longrightarrow}\quad(X(w),B).$$
The convergence in law is understood on $\mathbb{R}\times\mathcal{C}([0,1])$ and $X(w)$ denotes a random variable with the same law as that of  $X$ under $m$, function of a Brownian motion  $W$ independent of $B$.}\\

\noindent{\bf Proof.}  a) If $X$ has the form
$$X=\exp\{i\int_0^1\xi.dB+\frac{1}{2}\int_0^1|\xi|^2ds\}$$
for an element $\xi\in L^2([0,1],\mathbb{R}^d)$, we have
$$T_n(X)=\exp\{i\int_0^1\xi^\ast_sM_n(s)dB_s+\frac{1}{2}\int_0^1|\xi|^2ds\}.$$ where $\xi^\ast_s$ denotes the transposed of $\xi_s$. If we put $Z^n_t=\int_0^t\xi^t_sM_n(s)dB_s$ 
$$\langle Z^n,Z^n\rangle _t=\int_0^t\xi^t_sM_n(s)M_n^\ast(s)\xi_sds=\int_0^t|\xi|^2(s)ds$$ is a continuous function. Now by theorem 4
$$\int_0^t\xi^\ast_sM_n(s)ds\rightarrow\int_0^t\xi^\ast_sds\int_0^1M_n(s)ds=0.$$ Since the functions $t\mapsto \int_0^t\xi^\ast_sM_n(s)ds$ are uniformly continuous ($M$ is bounded) by Ascoli theorem $sup_t|\int_0^t\xi^\ast_sM_n(s)ds|\rightarrow0.$ The argument of Rootzen applies once more
$$(\int_0^.\xi^\ast M_ndB,B)\quad\stackrel{d}{\Longrightarrow}\quad(\int_0^.\xi.dW,B)$$ giving the result by the continuity of the exponential function.

b) For $X\in L^1(m)$, we consider $X_k$ linear combination of exponentials of the above form approximating $X$ in $L^1(m)$.

By a) we have $\forall h\in L^2([0,1],\mathbb{R}^d)$
$$\mathbb{E}[e^{iuT_n(X_k)}e^{i\int h.dB}]\rightarrow\mathbb{E}[e^{iuX_k}]\mathbb{E}[e^{i\int h.dB}]$$
but
$$|\mathbb{E}[e^{iuT_n(X)}e^{i\int h.dB}-\mathbb{E}[e^{iuT_n(X_k)}e^{i\int h.dB}]|
\leq |u|\mathbb{E}|T_n(X)-T_n(X_k)|=|u|\;\|X-X_k\|_{L^1}$$
what gives the result.

c) This extends to $\tilde{m}\ll m$ by the properties of stable convergence\hfill{$\Box$}\\

\noindent{\bf III.5. Stochastic differential equations from dynamics.}

In the case $f(x)=\theta(x)=\frac{1}{2}-\{x\}$, the approximation used in parts III.1 and III.2  approaches $B_t$ by 
\begin{equation}
B_t-\int_0^t(s-\frac{1}{2n}-\frac{[ns]}{n})dB_s
\end{equation}
and yields limit of the type
$$
(n(B_.-B_.^n),B_.)=(n\int_0^.(s-\frac{1}{2n}-\frac{[ns]}{n})dB_s,B_.)
\stackrel{d}{\Longrightarrow}(\frac{1}{\sqrt{12}}W_.,B_.)
$$
and
\begin{equation}
\left(n\int_0^.(s-\frac{[ns]}{s})dB_s,n\int_0^.(B_s-B_{\frac{[ns]}{n}})ds,B_.\right)\stackrel{d}{\Longrightarrow}(\frac{1}{\sqrt{12}}W_.+\frac{1}{2}B_.,-\frac{1}{\sqrt{12}}W_.+\frac{1}{2}B_.,B_.)
\end{equation}
Now, when we solve   by the Euler method a stochastic differential equation of the type defining a diffusion process and expand the coefficients in series, we encounter integrals of the type
$$\int_0^.(s-\frac{[ns]}{n})dB_s,\quad\int_0^.(B_s-B_\frac{[ns]}{n})ds$$
but also of the type
\begin{equation}
\int_0^.(B_s-B_{\frac{[ns]}{n}})dB_s
\end{equation}
and these last ones, by the central limit theorem, yield the convergence
\begin{equation}
\left(\sqrt{n}\int_0^.(B_s-B_{\frac{[ns]}{n}})dB_s,B_.\right)\stackrel{d}{\Longrightarrow}(\frac{1}{\sqrt{2}}\widetilde{W_.},,B_.).
\end{equation} 
Then let us remark that

a) the limits (18) are generally hidden by the limits (20) because of the order of magnitude of the coefficients $n$ and $\sqrt{n}$ respectively,

b) in (20), the conditional law of the random variable $\int_0^.B_{\frac{[ns]}{n}}dB_s$ with values in $\mathcal{C}([0,1],\mathbb{R})$ given $\int_0^.B_sdB_s$ is not reduced to a Dirac mass, the approximation is not deterministic (as can be seen, for instance, by changing the sign of a Brownian path after the time $T=\inf\{s\rangle \frac{n-1}{n}, B_s=0\}$, the $\sigma$-field generated by $\int_0^.B_sdB_s$ being $\sigma(B_s^2, s\leq 1)$ ). In (17) instead, $B^n$ is a deterministic function of $B$.\\

Nevertheless, for some stochastic differential equations the limits (18) remain dominant and determine the convergence. This concerns, for instance, stochastic differential equations of the form
\begin{equation}\left\{
\begin{array}{l}
X^1_t=x^1_0+\int_0^tf^{11}(X^2_s)dB_s+\int_0^tf^{12}(X^1_s,X^2_s)ds\\
X^2_t=x^2_0+\int_0^tf^{22}(X^1_s,X^2_s)ds
\end{array}\right.
\end{equation}
where $X^1$ is with values in $\mathbb{R}^{k_1}$,  $X^2$ in $\mathbb{R}^{k_2}$,  $B$ in $\mathbb{R}^d$ and $f^{ij}$ are matrices with suitable dimensions. Such equations are encountered  to describe the movement of mechanical systems under the action of forces with a random noise, when the noisy forces depend on the position of the system and the time. Typically 
$$\left\{
\begin{array}{l}
X_t=X_0+\int_0^tV_sds\\
V_t=V_0+\int_0^ta(X_s,V_s,s)ds+\int_0^tb(X_s,s)dB_s
\end{array}\right.
$$
which is a perturbation of the equation $\frac{d^2x}{dt^2}=a(x,\frac{dx}{dt},t)$. In such equations the stochastic integral may be understood as Ito as well as Stratonovitch. For the equation (21) the iterative method of Kurtz and Protter [14] (see also Jacod and Protter [11]) may be applied starting with the results obtained in generalizing of the arbitrary functions principle. This  yields the following result that we state in the case $k_1=k_2=d=1$ for simplicity.\\

\noindent{\bf Theorem 11.} {\it If functions $f^{ij}$ are $\mathcal{C}^1_b$, and if $X^n$ is the solution of {\rm(22)} by the Euler scheme, 
$$(n(X^n-X),X,B)\stackrel{d}{\Longrightarrow}(U,X,B)
$$
where the process $U$ is solution of the stochastic differential equation
$$
U"(t)=\sum_{k,j}\int_0^t\frac{\partial f^{ij}}{\partial x_k}(X_s)U^k_sdY^j_s-\sum_{k,j}\int_0^t\frac{\partial f^{ij}}{\partial x_k}(X_s)\sum_m f^{km}(X_s)dZ^{mj}_s
$$ where $Y_s=(B_s,s)^t $ and 
$$\begin{array}{l}
dZ^{12}_s=\frac{1}{\sqrt{12}}dW_s+\frac{1}{2}dB_s\\
dZ^{21}_s=-\frac{1}{\sqrt{12}}dW_s+\frac{1}{2}dB_s\\
dZ^{22}_s=\frac{ds}{2}
\end{array}
$$  
and as ever $W$ is an independent Brownian motion.}\\

Thus the Euler scheme for solving this kind of equations encountered in mechanics gives rise to an asymptotic weak limit, but  in $\frac{1}{n}$ and based on the arbitrary functions principle, instead of being in $\frac{1}{\sqrt{n}}$ and based on a version of the central limit theorem.

\begin{list}{}
{\setlength{\itemsep}{0cm}\setlength{\leftmargin}{0.5cm}\setlength{\parsep}{0cm}\setlength{\listparindent}{-0.5cm}}
  \item\begin{center}
{\small REFERENCES}
\end{center}\vspace{0.4cm}
 [1] {\sc Borel E.} {\it Calcul des probabilit\'es} Paris, 1924.

 [2] {\sc Bouleau N.} {\it Error Calculus for Finance and Physics, the Language of Dirichlet Forms}, De Gruyter, 2003.

[3] {\sc Bouleau N.}  ``When and how an error yields a Dirichlet form" {\it J. Funct. Anal.}, Vol 240, n$^0$2, (2006), 445-494.

[4] {\sc Fr\'echet M.} ``Remarque sur les probabilit\'es continues"{\it Bull. Sci. Math.} $2^e$ s\'erie, 45, (1921), 87-88.

[5] {\sc Fukushima, M., Oshima, Y., Takeda, M.} {\it Dirichlet forms and symmetric Markov processes}, De Gruyter 1994.

 [6] {\sc Hopf E.} ``On causality, statistics and probability" {\it J. of Math. and Physics} 18 (1934) 51-102.

[7] {\sc Hopf E.} ``\"{U}ber die Bedeutung der willk\"{u}rlichen Funktionen f\"{u}r die Wahrscheinlichkeitstheorie" {\it Jahresbericht der Deutschen Math. Vereinigung} XLVI, I, 9/12, 179-194, (1936).

 [8] {\sc Hopf E.} ``Ein Verteilungsproblem bei dissipativen dynamischen Systemen" {\it Math. Ann.} 114, (1937), 161-186.

 [9] {\sc Hostinsk\'y B.} ``Sur la m\'ethode des fonctions arbitraires dans le calcul des probabilit\'es" {\it Acta Math.} 49, (1926), 95-113.

 [10] {\sc Hostinsk\'y B.} {\it M\'ethodes g\'en\'erales de Calcul des Probabilit\'es}, Gauthier-Villars 1931.

 [11] {\sc Jacod, J., Protter, Ph.} ``Asymptotic error distributions for the Euler method for stochastic 
differential equations'' {\it Ann. Probab.} 26, 267-307, (1998)

 [12] {\sc von Kries, J.} {\it Die Prinzipien der Wahrscheinlichkeitsrechnung}, Freiburg 1886.

[13] {\sc Kahane j.-P., Salem R.} {\it Ensembles parfaits et s\'eries trigonom\'etriques}, Hermann (1963).

 [14] {\sc Kurtz, Th; Protter, Ph.} ``Wong-Zakai corrections, random evolutions and simulation schemes for SDEs" {\it Stochastic Analysis} 331-346,
Acad. Press, 1991.

 [15] {\sc Lyons, R.}``Seventy years of Rajchman measures" {\it J. Fourier Anal. Appl.} Kahane special issue (1995), 363-377.

 [16] {\sc von Plato, J.} ``The Method of Arbitrary Functions" {\it Brit. J. Phil. Sci.} 34, (1983), 37-42.

[17] {\sc Poincar\'e, H.} {\it Calcul des Probabilit\'es} Gauthier-Villars, 1912.

[18] {\sc Rajchman A.} ``Sur une classe de fonctions \`a variation born\'ee" {\it C. R. Acad. Sci. Paris}187, (1928), 1026-1028.

[19] {\sc Rajchman A.}``Une classe de s\'eries g\'eom\'etriques qui convergent presque partout vers z\'ero" {\it Math. Ann.}101, (1929), 686-700.

[20] {\sc Revuz D., Yor M.} {\it Continuous Martingales and Brownian motion}, Springer 1994.

 [21] {\sc Rootz\'en, H.} ``Limit distribution for the error in approximation of stochastic integrals'' {\it Ann. Prob. } 8, 241-251, (1980).

\end{list}

\end{document}